\documentclass[10pt]{article}
\usepackage{amssymb}
\begin{document}
\title{Finite groups acting on coherent sheaves\\ and Galois covers}
\author{Armando S\'{a}nchez-Arg\'{a}ez
\thanks{The author was partially supported by CONACyT project 27062-E.\newline
} }
\date{}
\maketitle

\newtheorem{Lema}{Lemma}
\newtheorem{Teo}{Theorem}

\newtheorem{Prop}{Proposition}
\newtheorem{Def}{Definition}
\newtheorem{Cor}{Corollary}

\newcommand{\mt}{\ensuremath{\mathcal{O}_{Y}[G]}}
\newcommand{\m}[1]{\ensuremath{\mathcal{#1}}}
\newcommand{\mo}{\ensuremath{\mathcal{O}}}
\newcommand{\me}{\ensuremath{\mathcal{E}}}
\newcommand{\mf}{\ensuremath{\mathcal{F}}}
\newcommand{\ml}{\ensuremath{\mathcal{L}}}
\newcommand{\ra}{\ensuremath{\longrightarrow}}
\newcommand{\da}{\ensuremath{\downarrow}}
\newcommand{\st}[1]{\ensuremath{\stackrel{#1}{\longrightarrow}}}
\newcommand{\wh}[1]{\ensuremath{\widehat{#1}}}
\newcommand{\id}[1]{\ensuremath{{#1}_{*}}}
\newcommand{\ii}[1]{\ensuremath{{#1}^{*}}}
\newcommand{\pit}{\ensuremath{{\pi}^{*}}}
\newcommand{\pib}{\ensuremath{{\pi}_{*}}}
\newcommand{\gt}{\ensuremath{{g}^{*}}} \newcommand{\gb}{\ensuremath{{g}_{*}}}
\newcommand{\tn}[2]{\ensuremath{{#1}\otimes {#2}}}
\newcommand{\oxg}{$\mathcal{O}_{X}[G]$-module}
\newcommand{\oxa}{$\mathcal{O}_{X}(A)$-module}
\newcommand{\oyg}{$\mathcal{O}_{Y}(G)$-module}
\newcommand{\ox}{$\mathcal{O}_{X}$-module}
\newcommand{\qa}{\'{a}}
\newcommand{\qe}{\'{e}}
\newcommand{\qi}{\'{\i}}
\newcommand{\qo}{\'{o}}
\newcommand{\qu}{\'{u}}
\newcommand{\mi}[1]{\ensuremath{\mathit{#1}}}

\begin{abstract}
Let $G$ be a finite group and
$\rho:G\ra End(\me)$ be a group representation of $G$ on a coherent sheaf over
an integral scheme. The purpose of this paper shall give a decomposition
theorem of such representations in non-splitting components and apply this
results
to the studie of Galois covers of a variety.
\end{abstract}

\section{Introduction}
In this paper by a variety we mean an integral separated scheme of finite type
over an
algebraically closed field $k$. All the varieties are projective. All the
groups that we will consider are finite and $char(k)\mid \hspace{-2.4mm}/
\mid G \mid$.

Let $X$ be a variety and $G$ be a finite group, let $p:E\ra X$ be a vector
bundle with $G$ as subgroup of $Aut_{X}(E)$. Then we have a representation
$p^{-1}(x)$ of $G$ over $k$ for any closed point $x\in X$. This is a natural
generalization of the representation theory of $G$ over $k$; in fact we
recover this last when $X=Spec(k)$. Now, to have a representation $G$ on the
vector bundle is equivalent to have a structure of \oxg \, in the sheaf of
sections of $E$. Then in this paper we focus our attention to the abelian
category of free torsion coherent sheaves with a group action. The main
purpose of the first three section shall give a classification theorem of such
objects.

The first motivation of that studie is the relationship between Galois
covers of varieties and representation theory over a field, for see this
let $Y$ be a projective variety and $G$ be a finite group acting on $Y$, let
$\pi : Y \ra Y/G=X $ the quotient map. Then $\pib \mo_{Y}$ is a coherent free
torsion $\mo_{X}$-module, such that, for any open set $U\subset X$, $\pib
\mo_{Y}(U)$ is an $\mo_{X}(U)[G]$-module, furthermore in the generic point
$\pib \mo_{ Y, \epsilon}$ is isomorphic to $K_{X}[G]$ as representation of $G$
over the rational function field $K_{X}$ of $X$. In section 5, structure
theorems for
$\pib \mo_{X}$ are given, many of them are generalizations from the cyclic
case.

\section{General Theorems.}

Let $A$ be a semisimple $k$-algebra  (not necessary commutative)
of finite dimension  over $k$. It is very known that
$A=e_{1}A\oplus ....\oplus e_{r}A$
with $e_{i}A$ a simple algebra and $e_{i}$ idempotent, let $n_{i}=dim_{k}
e_{i}A$.

\begin{Def}
Let $(X,\m{O}_{X})$ be a ringed space. Define $\mathcal{O}_{X}(A):=
\mathcal{O}_{X}\otimes_{k}A$, this is a sheaf of $\mo_{Y}$-algebras, not
necessary
commutative, over $X$ but as an  $\mathcal{O}_{X}$-module is free of rank
$n=dim_{k}A$. An  \oxa\ is a pair consisting of an \ox\
$\m{E}$ together with a $k$-morphism of rings $\rho:A \rightarrow
End(\m{E})$. A morphism of   \oxa s is a morphism
$\phi:\m{E}\rightarrow \m{F}$ of \ox s such that the next diagram commute
for any $a\in A$ \[ \begin{array}{rcccl}\; & \m{E} & \stackrel{\phi}
\rightarrow & \mathcal{F} & \; \\ \rho (a) &\downarrow & \; & \downarrow &
\rho (a)\\ \; & \m{E} & \stackrel{\phi}\rightarrow & \mathcal{F} &
\;\end{array} \] In this case we say that $\phi$ is $A$-invariant
\end{Def}

\vspace{3mm}

\textbf{Remark and notation 1} Clearly from this definition,
the kernel, cokernel, and image of a morphism of
$\mo_{X}(A)$-modules is again  an
$\mo_{X}(A)$-module. Any direct sum, direct product, direct limit or inverse
limit of $\mo_{X}(A)$-modules is an $\mo_{X}(A)$-module. If $x\in X$,
$\m{E}_{x}$ have a natural structure of $\m{O}_{X,x}(A)$-module.
Observe,that this is an abelian category.

If $\m{E}$ and $\m{F}$ are
\oxa s, we denote by $Hom_{A}(\m{E},\m{F})$ the group formed by
the morphisms of \oxa s. Let $U
\subset X$ be an open set of $X$ and $\m{E} ,\m{F}$ two \oxa s then
 we define by
 $\m{H}\mathit{om}_{A}(\m{E},\m{F})$ the sheaf  $$U\mapsto
Hom_{A}(\m{E}_{|U},\m{F}_{|U})$$
although, it is a sheaf of \ox s, it is not necessarily of
 \oxa s.

\vspace{3mm}

Recall that $A$ is a semisimple algebra, then we have only a finite number of
simple $A$-modules, let
$V_{1},...V_{r}$ be this modules. Now is very known that
for any finitely generated $A$-module $M$ it decompose in a direct
sum
 $M_{1}\oplus...\oplus M_{r}$ where each
$M_{i}$ is isomorphic to $V_{i}^{n_{i}}$ and there exist
$\{e_{1},...e_{r}\}\subset A$ such that $e_{i}e_{j}=e_{i}\delta _{i,j},
e_{1}+...+e_{r}=1$ and $e_{i}M=M_{i}$.

\begin{Prop}
If $\mathcal{M}$ is an $\mathcal{O}_{X}(A)$-module, not necessarily
finite generated, then $\mathcal{M}= \mathcal{M}_{1}\oplus ....\oplus
\mathcal{M}_{r}$ where $\mathcal{M}_{i}=e_{i} \mathcal{M}$ , furthermore if
$\mathcal{M}$ is locally free then each $\mathcal{M}_{i}$
is locally free.$\diamondsuit$
\end{Prop}

\begin{Def} Let $\m{E}$ an \oxa , we define the isotypical
decomposition
of
 $\m{E}$,
by the decomposition
$$\m{E}=\m{E}_{1}\oplus...\oplus \m{E}_{r}$$
obtained in the above proposition.
\end{Def}

\begin{Prop} Let $\m{E}$ and $\m{F}$ be two \oxa s and $\phi :\m{E}
\rightarrow \m{F}$ be an $\m{O}_{X}(A)$-morphism. If
$e_{1}\m{E}\oplus...\oplus e_{r}\m{E}$ and $e_{1}\m{F}\oplus...\oplus
e_{r}\m{F}$
are the isotypical decomposition of
$\m{E}$ and $\m{F}$ respectively,
then $\phi$ decompose in morphisms
$$\phi_{i}:e_{i}\m{E} \rightarrow e_{i}\m{F}$$ $\diamondsuit$
\end{Prop}

\begin{Cor} Let $\m{E}$ and $\m{F}$ be two \oxa s and $\phi :\m{E}
\rightarrow \m{F}$ be an $\m{O}_{X}(A)$-morphism. Then $\m{H}om_{A} (\m{E},
\m{F})=\bigoplus \m{H}om_{A} (e_{i}\m{E}, e_{i}\m{F})$.
$\diamondsuit $ \end{Cor}

\begin{Cor} Any exact sequence of \oxa s
$$...\rightarrow \m{E}_{i-1}\rightarrow \m{E}_{i}\rightarrow
\m{E}_{i+1}\rightarrow ...,$$
 decompose in exact sequences
$$...\rightarrow e_{j}\m{E}_{i-1}\rightarrow
e_{j}\m{E}_{i}\rightarrow e_{j}\m{E}_{i+1}\rightarrow ...$$
with $j\in \{1,...,r\}. \; \diamondsuit$
\end{Cor}

\section{$\mathcal{O}_{X}[G]$-Modules.}

Since we are interested on sheaves with group action,
from now on we suppose that $A$ is the group algebra
$k[G]$, where
$char(k)\mid \hspace{-2.4mm}/ |G|$.
We denote by
$\mathcal{O}_{X}[G]:=\mathcal{O}_{X}(k[G])$ and by $V_{0},...,V_{r}$ the
irreducibles representations of $G$ over $k$, this are the simples
$k[G]$-modules, and by $\mo(V)$ the $\mo_{X}[G]$-module definite
by $\mo_{X}\otimes_{k}V$, where $V$ is a representation of $G$ over
$k$. We will use $V_{0}$  to mean the trivial representation of
$G$.

\vspace{3mm}

Under the hypothesis $A=k[G]$, if $\m{E}$ and $\mf$ are \oxg s, then
the $\mathcal{O}_{X}$-modules
$\mathcal{E} \otimes \mathcal{F}$ and $\mathcal{H}\mathit{om}(\mathcal{E} ,
\mathcal{F})$ have a natural structure of
$\mathcal{O}_{X}[G]$-modules, given by

$$g(x\otimes y)=(gx)\otimes (gy),\; (\forall \; g\in G, x \in \m{E}, y\in
\m{F})$$
and
 $$(g\sigma )(e)=g(\sigma )(g^{-1}e)$$
respectively

\vspace{3mm}

\textbf{Notation 2:} By convention, a locally free \oxg\ $\me$ means that it is
locally free as an \ox .

\begin{Lema}\label{identi}
Let $\mathcal{E}$ be a locally free $\mathcal{O}_{X}[G]$-module of finite rank.
Give us to $\m{E}^{\vee}$ the natural structure of \oxg . Then
we have the next natural isomorphisms between
\oxg s:

a)$(\mathcal{E}^{\vee})^{\vee} \cong \mathcal{E}$.

b)For any $\mathcal{O}_{X}[G]$-module $\mathcal{F}$,
$\mathcal{H}\mathit{om} (\mathcal{E,F})\cong \mathcal{E}^{\vee }\otimes
\mathcal{F}$

c)For any $\mathcal{O}_{X}[G]$-modules $\mathcal{F,G}$,
$Hom_{\mathcal{O}_{X}}(\mathcal{E}\otimes
\mathcal{F,G})$  $\cong$
$Hom_{\mathcal{O}
_{X}}(\mathcal{F},\mathcal{H}\mathit{om} (\mathcal{E,G}))$.$\diamondsuit$
\end{Lema}

\vspace{3mm}

Now we procede to proof the main theorem of this section,

\begin{Teo}\label{descom} Let $X$ be an integral scheme over an algebraically
closed field $k$, let
$K$ be the rational function field of $X$ and  $\epsilon$ the generic point.
Let $\m{E}$ be a free torsion coherent \oxg ,
If $\m{E}_{\epsilon}\simeq (V_{0}^{n_{0}}\oplus....\oplus V_{r}^{n_{r}})\otimes
_{k}K$ is the representation of $G$  in the generic point of $X$, then,
the natural decomposition $$\m{E}=e_{0}\m{E}\oplus....\oplus e_{r}\m{E}$$
satisfies

i) $rank(e_{i}\m{E})=n_{i}\times dim(V_{i})$

ii) The natural inclusion $(e_{i}\m{E})\hookrightarrow
\m{E}$ induce an isomorphism of representations
$(e_{i}\m{E})_{\epsilon}\hookrightarrow V_{i}^{n_{i}}\otimes _{k}K$
\end{Teo}
\textit{Proof:}
As the question is local, let be $U=Spec(A)$ an affin open set of $X$ where
$\m{E}$ is trivial
 i.e
$\mathcal{E}_{|_{U}}\cong A^{n}$,
and consider the isotypical decomposition
$$\hspace{-45mm}(1)\hspace{40mm}\mathcal{E}=\mathcal{E}_{1}\oplus ....\oplus
\mathcal{E}_{r}$$
where $\mathcal{E}_{i|_{U}}\cong
e_{i}\mathcal{E}{|_{U}}$, on the other hand $X$ in an integral scheme, and so
 $A$ is an integral domain and $K(X)$ is the quotient field of $A$, taken 
localization in
the generic point
we have that $\mathcal{E}_{X}=A^{n}\otimes
_{A}K(X)=K(X)^{n}$. Now, the next commutative diagram

\begin{picture}(150,70)
\put(105,45){$G$}
\put(145,45){$Aut(K(X)^{\oplus n})$}
\put(145,5){$Aut(A^{\oplus n})$}
\put(165,18){\vector(0,1){20}}
\put(120,40){\vector(1,-1){25}}
\put(120,50){\vector(1,0){20}}
\end{picture}

\hspace{-6mm}say us that the decomposition (1) is determinate in the generic
point by the representation of
$G$ in $K(X)^{\oplus n}$, from this, the decomposition of
 $\mathcal{E}$ is determinate by the given in the generic point as
representation of
$G$. Further, the action of
$G$ on $B^{n}$ is given by the restriction of the action of
 $G$ on $K(X)^{n}$.
 Now, the theorem follows from the next
\begin{Lema} Let $k$ be an algebraically closed field and $G$ be a finite group with
$char(k) | \hspace{-1.25mm}/ |G|$. Then, for any field extention $K$ of $k$, a
representation of $G$ over $K$ is of the form $V\otimes_k K$ with $V$ a 
representation of
$G$ over $k$.
\end{Lema}
\textit{Proof:} It is sificient to prove the lemma for irreducible representations. 
By the
corollary 3.61 in \cite{cur} page 68, any  representation of the form $V\otimes_k K$ 
is
irreducible over $K$ if $V$ is irreducible over $k$, and by theorem 30.15 in 
\cite{cur2}
page 214 all this are different irreducible representations, then this are all the
irreducible representations.
$\diamondsuit$

\begin{Def}Let $X$ be an integral scheme over an algebraically
closed field $k$,
$K$ be the function field of $X$ and  $\epsilon$ be the generic point.
Let $\m{E}$ be a free torsion coherent \oxg ,
and
$\m{E}_{\epsilon}\simeq V\otimes_{k}K$ be the representation of $G$ in the
generic point, we define the type of the representation of $G$ on $\me$ as the
representation $V$
\end{Def}

\begin{Cor} Let $X$ and $\me$ be as in the above theorem and $V$ be the
type of the representation in $\me$. Then the representation of $G$ in
the fiber of a closed point of $X$ is generically $V$.
\end{Cor}

\begin{Cor} Let $X$ be an integral scheme over an algebraically closed field
$k$, If $\m{E}$ and $\m{F}$ are \oxg s with type $V$ and $W$ respectively,
Then

$i)\;\m{E} \otimes
_{\m{O}_{X}}\m{F}$ have the type $V\otimes W$.

$ii)\;\Lambda^{r}\m{E}$ have the type $\Lambda ^{r}V$, and

$iii)\; S^{n}\m{E}$ have the type $S^{n}V$
\end{Cor}
\textit{Proof:}
Just we need to observe that if $\epsilon$ is the generic point of $X$
then $(\me \otimes_{\mo_{X}} \mf)_{\epsilon}=
\me_{\epsilon} \otimes_{K(X)} \mf_{\epsilon}$, $(\Lambda^{r}\m{E})_{\epsilon}
= \Lambda^{r}(\m{E}_{\epsilon})$ and $ (S^{n}\m{E})_{\epsilon}=
S^{n}(\m{E}_{\epsilon})$, the conclusion is immediately from the general theory
of representation over a field.$\diamondsuit$

\begin{Cor}
Let $\mathcal{E}_{i},\mathcal{E}_{j}$ be two locally
free $\mathcal{O}_{X}[G]$-modules  of type $V_{i},V_{j}$, respectively,
both different irreducible representations of $G$. Then
$\mathcal{H}\mathit{om}_{G}(\mathcal{E}_{i},\mathcal{E}_{j})=0$
\end{Cor}
\textit{Proof:}
Let $\epsilon$ be the generic point of $X$. From the theorem 1, the rank of
$\mathcal{H}\mathit{om}_{G}(\mathcal{E}_{i},\mathcal{E}_{j})$ is the dimension
of the part of type $V_{0}$ of
$(\mathcal{E}_{i}^{\vee}\otimes
_{\mathcal{O}_{X}} \mathcal{E}_{j})_{\epsilon}$, but it is zero by the general
theory of representation over a field. Now
$\mathcal{H}\mathit{om}_{G}(\mathcal{E}_{i},\mathcal{E}_{j})$ is a locally
free sheaf, then we conclude.
$\diamondsuit$

\begin{Cor}
Let $k=\mathbb C $ and $\mathcal{E}_{1},\mathcal{E}_{2}$ be two locally free
$\mathcal{O}_{X}[G]$-modules with type $V_{i}$. Then

$$rank \mathcal{H}\mathit{om}_{G}(\mathcal{E}_{1},\mathcal{E}_{2})=
rank\ \mathcal{E}_{1}\times  rank\ \mathcal{E}_{2}/dimV_{i}^{2}$$

\end{Cor}

\section{Irreducibles $\mo_{X}[G]$-Modules }

Let $\mf$  be a free torsion sheaf over an integral scheme. A natural question
is: when $\mf$ admit an
$\mo_{X}[G]$-module structure? As we see posterior, this will be possible
only if the decomposition of  $\mf$ have certain structure.

\begin{Lema} Let $V$ be an irreducible representation of $G$ over $k$, and
suppose that $dim_{k}V\geq 2$. Then, there are a subgroup
 $H \leq G$ such that $V\da _{H}^{G}$ have at least two
isotypical
components
  \end{Lema}
\textit{Proof:} Suppose that it is false. Then
for any $g\in G$, \, $V\da _{<g>}$ have only one isotypical
component,
that
means,  any element $g$ acts on $V$ by multiplication of
a constant, in this case, any subspace $V$ of dimension 1 is $G$ invariant,
but this is a contradiction with our hypothesis because  $V$ is irreducible
of dimension great than 1.$\diamondsuit$

\begin{Lema}  Let $\mf$ be a non-splitting coherent \ox\ , if $\mf$ is an
\oxg\ for some finite group $G$, then the representation type is
$V=W^{n}$ for some irreducible $W$ of dimension 1. Furthermore, any
 element $g\in G$ acts by constant multiplication .
\end{Lema} \textit{Proof:}
As $\mf$ is irreducible, the type of representation must be
 $V=W^{n}$ for some irreducible
  $W$. Furthermore, when we restrict the action to any subgroup of
 $G$, the structure must be preserved. But, by the above lemma
 $W$ must have dimension 1, and in particular, if
$g\in G$,   this acts in the generic point by a constant multiplication, then
acts globally in the same way.
$\diamondsuit$

\vspace{3mm}

\textbf{Remark 3} Let $V$ be an irreducible representation of $G$ over $k$. Then by 
general theory, de dimension of the type $V_0$ part of $V^s\otimes V^*$ is equal to 
$s$, where $V^*$ is the dual representation of $V$. Thus by lemma 2 the same is 
true for representation over extensions of $k$.

\textbf{Remark 4}Using the lemma 2, we obtain that the number of irreducible 
representations of $G$ over $K$ is the same than over $k$.

Now we are in position
of classify the non-splitting
$\mathcal{O}_{X}[G]$-modules

\begin{Teo}
Let $X$ be an integral scheme, let
 $\mathcal{E}$ be a non-splitting free torsion coherent
$\mathcal{O}_{X}[G]$-module
of $W$ type. Then $\mathcal{E} \simeq
\mathcal{O}(V)\otimes \mathcal{F}$ with $\mathcal{F}$ a non-splitting
$\mathcal{O}_{X}$-module and $W \simeq V^{rank \mf}$, with $V$ an
irreducible representation. \end{Teo}
\textit{Proof}:
Let $\m{E}=\m{F}_{1}^{n_{1}}\oplus...\oplus\m{F}_{r}^{n_{r}}$ be
the decomposition of
 $\m{E}$ on non-splitting
 $\m{O}_{X}$-modules with
$\m{F}_{i}\neq\m{F}_{j}$, if $i\neq j$. This decomposition is unique up to
permutations. From this, if
 $g\in G$,
$g\m{E}=g\m{F}_{1}^{n_{1}}\oplus...\oplus g\m{F}_{r}^{n_{r}}$, then
$g\m{F}_{i}=\m{F}_{j}$  for some
$j$, but this imply
 $\m{F}_{i}\simeq \m{F}_{j}$ and
$i=j$, Then, any
$\m{F}_{i}^{n_{i}}$ is
$G$-invariant but
$\me$ is a non-splitting \oxg , and so
 $r=1$, i.e.
$\m{E}=\m{F}^{n}$ with
 $\m{F}$ a non-splitting
$\m{O}_{X}$-module. Therefore the type of the representation is
$V^{s}$ with
$V$ irreducible.

The next step is to show that
$s=rank \, \mf$ and $dim \, V=n$.

Now, we consider the part of type
 $V_{0}$ of
 $\mo(V^{\vee})\otimes \me$, where
 $V^{\vee}$ is the dual representation of
 $V$, This is a direct component of
 $\mo(V^{\vee})\otimes \me\simeq \mf^{ndimV}$, so this component must be
 $\mf^i$ for some $i$.
In other hand, let $\epsilon$ be the generic point of $X$, so $\me_{\epsilon}=
V^{s}\otimes K(X)$, then by remark 3, $s=irank \m F$. Now, let
$g\in G$, and consider the representation of
 $<g>$ obtained by restriction, then
 $\me$ have an isotypical decomposition
 $\m{E}=\m{F}_{\chi_{1}}^{n_{1}}\oplus...\oplus\m{F}_{\chi_{r}}^{n_{r}}$
 with $\chi_{i}$ the irreducible representations of
$<g>$. Then, by the above lemma, in each component
$<g>$ acts by constant multiplication, then
$\me=[\mo(V_{\chi_{1}})^{n_{1}}\oplus...\oplus\mo(V_{\chi_{r}})^{n_{r}}]\otimes
\mf $, and in the generic point we must have
$\me_{\epsilon}=
[V_{\chi_{1}}^{n_{1}}\oplus...\oplus V_{\chi_{r}}^{n_{r}}]\otimes
\mf_{\epsilon}$, in other hand, the restricted representation is given by
$\me_{\epsilon}=(V\downarrow ^{G}_{<g>})\otimes K(X)^s$, now using the
uniqueness of the decomposition in irreducible representations over a field,  we 
obtain
 $(V\downarrow ^{G}_{<g>})^i\otimes K(X)^{rank \m F}=
[V_{\chi_{1}}^{n_{1}}\oplus...\oplus V_{\chi_{r}}^{n_{r}}]\otimes
\mf_{\epsilon}$
and using
the fact that the characters are a basis for
$K[G]$ when $K$ is an extension of an algebraically closed field
(see remark 4 and \cite{cur2} page 213 Theorem 30.12), then 
$[V_{\chi_{1}}^{n_{1}}\oplus...\oplus V_{\chi_{r}}^{n_{r}}]\simeq (V\downarrow ^G 
_{<g>})^i \simeq (V^i)\downarrow ^G _{<g>}$. Then $\m E =\mo (V^i)\otimes \m F$, but 
$\m E$ is a non-splitting $\mo_X(G)$-module, so $i=1$ and in consequence,
$s=rank\mf , \,dimV=n$ and $\me\simeq \mf^{dimV}$.

Now, in the last step we want to describe the representation of $G$ on $\me$.
For this, we consider again the part of type
 $V_{0}$ of
 $\mo(V^{\vee})\otimes \me$, where
 $V^{\vee}$ is the dual representation of
 $V$. This is a direct component of
 $\mo(V^{\vee})\otimes \me\simeq \mf^{n^{2}}$ and by the above considerations, his 
rank is
 $rank \, \mf$,and so this component must be
 $\mf$. Consider the natural inclusion
$$ \mo(V)\otimes\mf \ra \mo(V)\otimes\mo(V^{\vee})\otimes\me$$
and the morphism
$$\mo(V)\otimes\mo(V^{\vee})\otimes\me\ra \me$$
given by
 $$a\otimes\delta\otimes e\mapsto \delta(a)e$$
then we have the natural $G$-morphisms
$$\mo(V)\otimes \mf \ra \me\simeq\mf^{dimV} $$
that is a $G$-isomorphisms in the generic point, and using the integral
hypothesis we obtain that the kernel is a torsion subsheaf of $\me$, but it is a free
torsion coherent sheaf by hypothesis, so we have that the morphism is injective. Now 
using
the hypothesis of projective, we can consider the Hilbert polinomial of the cokernel 
(with
respec to a some ample sheaf), then it must be zero, so the morphism is surjective 
and
then it is a $G$-isomorphism. $\diamondsuit$

\vspace{3mm}

Now we have the classification theorem of \oxg .
\begin{Teo}
Let $\m{E}$ be a free torsion coherent  \oxg\  of type
$V_{0}^{n_{0}}\oplus...\oplus
V_{r}^{n_{r}}$. Then the isotypical decomposition of $\m{E}$ is
given
by
$$\m{E} \simeq \m{O}(V_{0})\otimes \m{F}_{0}\oplus...\oplus
\m{O}(V_{r})\otimes \m{F}_{r}$$
when $\m{F}_{i}$ is an \ox\ of rank $n_{i}\; \diamondsuit$
\end{Teo}

\begin{Prop} Let
$$\me_{i}=\mo(V_{0})\otimes \mf_{0,i}\oplus...\oplus \mo(V_{r})\otimes
\mf_{r,i},\, \, \, \, i=1,2.$$
Then
$$\m{H} om_{G}(\me_{1},\me_{2})\st{\sim}\bigoplus_{j=0}^{r}
\m{H}om_{\mo_{X}}(\mf_{j,1},\mf_{j,2})$$
\end{Prop}
\textit{Proof:}
By corollary 1
 $$Hom_{G}(\me_{1},\me_{2})= \bigoplus_{j=0}^{r}
Hom_{G} (\mo(V_{j})\otimes \mf_{j,1},\mo(V_{j})\otimes \mf_{j,2})$$
then we can suppose, without lost of generality, that
$\me_{i}=\mo(V)\otimes \mf_{i}$ with
$V$ an irreducible representation of $G$.

Let define the natural map
 $$\sigma :\mf_{1}\ra\mf_{2} \mapsto
Id\otimes \sigma:\mo(V)\otimes \mf_{1}\ra  \mo(V)\otimes \mf_{2}$$
and define the inverse map in the next way: Let
 $V^{\vee}$ be the dual representation of
$V$, then
 $V\otimes V^{\vee}\sim V_{0}\oplus W$, then
 $\mo(V)\otimes \mo(V^{\vee})\sim \mo_{X} \oplus \mo(W)$. Thus if
$\rho:
\mo(V)\otimes \mf_{1}\ra  \mo(V)\otimes \mf_{2}$ is a
$G$-morphisms we have a natural $G$-morphisms
$$Id\otimes
\rho:\mo(V^{\vee})\otimes \mo(V)\otimes \mf_{1}\ra \mo(V^{\vee})\otimes
\mo(V)\otimes \mf_{2}$$
let define $\psi$ by taken the restriction to the component corresponding to
the trivial type. It is clear that $\sigma$ and $\psi$
are inverse maps
.$\diamondsuit$

\begin{Prop}\label{inescin} Let $\me ,\mf_{0},...\mf_{r}$ be free torsion
coherent
sheaves over $X$, $W=V_{0}\oplus...\oplus V_{r}$ be
a $G$ representation over $k$ and suppose that
$$\phi :\mo(W)\otimes \me \ra  (\mo(V_{0})\otimes \mf_{0})\oplus...
\oplus(\mo(V_{r})\otimes \mf_{r})$$
is a $G$-invariant morphism. Then there exist a natural morphisms
$$\bar{\phi}:\me \ra
\mf_{0}\oplus...\oplus\mf_{r}$$
Moreover
 $\phi$ is injective if and only if $\bar{\phi}$ is.
\end{Prop}
\textit{Proof:}
By the above proposition
 $Hom_{G}(\mo(W)\otimes \me ,
(\mo(V_{0})\otimes \mf_{0})\oplus... \oplus(\mo(V_{r})\otimes \mf_{r}))$

is natural isomorphic to $\bigoplus_{i=0}^{r}Hom(\me, \mf_{i})$,
let $\bar{\phi}$ be the image of
 $\phi$ by this map .

Now
 $\phi$ is injective iff
$$\phi_{i}:\mo(V_{i})\otimes\me \ra \mo(V_{i})\otimes \mf_{i}$$ it is for
$i=0,...,r$, iff
$$Id\otimes \phi_{i}:\mo(V_{i}^{\vee})\otimes\mo(V_{i})\otimes\me \ra
\mo(V_{i}^{\vee})\otimes \mo(V_{i})\otimes \mf_{i}$$ it is for
 $i=0,...,r$, iff
 $$\bar{\phi_{i}}:\me\ra\mf_{i}$$ for
$i=0,...,r$ $\diamondsuit$

\section{Structure theorems for Galois covers.}
Let $X\ra Y$ be a Galois cover with $\m{G}al(X/Y)=G$. Then $\pib \mo_{X}$ have
a
natural structure of $\mo_{Y}[G]$-module, the objective of this section shall
give structure theorems for this object.

\begin{Teo}\label{descom} Let $X$ be a variety over an algebraically closed
field $k$
and suppose
$X$ with an action of a finite group $G$
($char(k)\mid \hspace{-2.8mm}/ \mid G\mid $)
Let $Y$ be the quotient variety. Then
 $$\pi_{*}\m{O}_{X}=\m{O}_{Y}\oplus (\m{O}(V_{1})\otimes
\m{E}_{V_{1}})\oplus ...\oplus (\m{O}(V_{r})\otimes \m{E}_{V_{r}}) $$
when $V_{i}$ are the irreducibles representation of $G$ over $k$ and $rank\,
\me_{V_{i}}=dim V_{i}$ \end{Teo}
\textit{Proof:}
Let $\pi:X\rightarrow X/G=Y$ be the natural quotient map, this is a Galois
cover of $Y$. Then the field extension $K(X):K(Y)$ is Galois with
Galois group $G$. Now, by the normal basis theorem, there exist an element
$\zeta \in K(X)$ such that $\{ g\zeta \} _{g\in G} $ is a basis of  $K(X)$ as
a vector space over $K(Y)$. Thus  $K(X)$ is the regular representation of $G$
over $K(Y)$, in other hand this is the action of $G$ on
$(\pi_{*}\m{O}_{X})_{\epsilon (Y)}$. Then the representation of $G$ in the
generic point of $Y$ is the regular one, and using the theorem 3 we have the
desired conclusion.$\diamondsuit$

\begin{Cor}\label{invar} Let $X$ and $Y$ be as above, let
 $H\leq G$ and $\phi: X/H \ra Y$ be the natural morphisms. Then
$$\phi_{*}\m{O}_{X/H}=(\pib \mo_{X})^{H}=\m{O}_{Y}\oplus
(\m{O}(V_{1}^{H})\otimes \m{E}_{V_{1}})\oplus ...\oplus
(\m{O}(V_{r}^{H})\otimes \m{E}_{V_{r}}). $$
\end{Cor}
\textit{Proof:} See \cite{mum} pages 65-69

\vspace{3mm}

Let $\pi: X\ra Y$ be a Galois cover with $\m{G}al(X/Y)=G$ and let
$\m{L}$ be an $\mo_{X}$-module ,
Let $\phi:\pib\m{L}\otimes \pib \mo_{X} \ra \oplus_{g\in G} \gt \pib
\ii{g}\m{L}$ be the map given
by
$$m\otimes b
\mapsto (m.g^{-1}b)_{g}.$$ Then we have

\begin{Lema} Let $X\ra Y$ be a Galois cover. Then
 $$0\ra\pib\pit\pib\mo_{X}\st{\phi_{*}} \bigoplus_{g\in G}
\pib (g^{*}\mo_{X}) \ra D\ra 0$$
is an exact sequence of
$\mo_{Y}(G)$-modules, where $supp D \subset
supp \pib\pit\pib\Omega_{X/Y}$.
\end{Lema}
\textit{Proof:} Just we need to proof the injective part of the sequence.

 Let $y\in Y-supp(\pib\Omega_{X/Y})$ a closed point of $Y$. Then the
localized sequence in $y$ is exact, so the sequence is exact and
$supp\, D \subset supp \,(\pib\Omega_{X/Y})$.$\diamondsuit$

\begin{Cor} Let
$\mf$ be a locally free sheaf  on $X$. Then
$$0\ra\pib\pit\pib\m{F}\ra \bigoplus_{g\in G} \pib (g^{*}\m{F}) \ra
D\ra 0$$
is an exact sequence of $\mo_{Y}(G)$-modules, where $supp\ D \subset
supp(\pib\Omega_{X/Y})$. \,$\diamondsuit$
\end{Cor}

Now we observe that $\pib\m{F}\otimes \pib \mo_{X}$ and $\bigoplus_{g\in
G}\pib\ii{g}\m{F}$
have a natural structure of \oyg \, given by
$$m\otimes a \stackrel{h}{\longmapsto} m\otimes h(a)$$
and
$$(m_{g})_{g}\stackrel{h}{\longmapsto} (m_{g})_{hg}$$
respectively, and the morphism $\phi$ is $G$-invariant.

\begin{Prop}\label{paso1} Let $X\st{\pi} Y$ be a Galois cover
with Galois group $G$, and $\mf$ be a locally free sheaf on $X$. Then we have
natural injective morphism
$$0\ra \pib
\m{F}\otimes_{\mo_{Y}}\me_{V_{i}}\st{\phi_{i}} \pib \m{F}^{dim V_{i}}.$$
\end{Prop}
\textit{Proof:}
Now by theorem \ref{descom}
$$\pi_{*}\m{O}_{X}=\m{O}_{Y}\oplus (\m{O}(V_{1})\otimes
\m{E}_{V_{1}})\oplus ...\oplus (\m{O}(V_{r})\otimes \m{E}_{V_{r}}) $$
is the isotypical decomposition of $\pib\mo_{X}$. Then we have that
$$\pib\m{F}\otimes \pib \mo_{X}=\pib\m{F}\oplus \bigoplus_{i=1}^{r}\otimes
\m{E}_{V_{r}}\otimes \m{F}$$ is the isotypical decomposition of it.

On the other hand, the action of $G$, in the generic point $\epsilon$ of $Y$ ,
in $\bigoplus_{g\in  G}\pib\ii{g}\m{F}$ is given by
$k[G]\otimes_{k}(\bigoplus_{g\in  G}\pib\ii{g}\m{F})_{\epsilon}$. Then the
isotypical decomposition of $\bigoplus_{g\in  G}\pib\ii{g}\m{F}$ is given by
$$\bigoplus_{g\in G}
\pib \gt \m{F}=\pib
\m{F}\oplus\mo
(V_{1})^{dimV_{1}}\otimes \pib
\m{F}\oplus
...\oplus\mo (V_{r})^{dimV_{r}}\otimes \pib
\m{F}$$
Now, recall that $\phi$ is $G$-invariant, then $\phi$ decompose on morphism
$$\phi:\mo(V_{i})\otimes[\me\otimes \pib\mf]\ra
\mo(V_{i})\otimes(\pib\mf)^{dim \, V_{i}}$$
and applying the proposition 4
we have the injective morphisms
$$\phi_{i}:\me\otimes \pib\mf\ra (\pib\mf)^{dim \, V_{i}}$$
and the proposition is proved. $\diamondsuit$

\vspace{3mm}

Now we are interested in the case when $\mf$ is $G$-stable (see \cite{mum}
pages 65-69). In this case the direct image of $\mf$  have a natural
structure of \oxg . So we have in $\pib\mf\otimes \pib\mo_{X}$ and
$\bigoplus_{g\in  G}\pib\ii{g}\m{F}$ two others actions of $G$ given by
$$m\otimes a \stackrel{\hat{h}}{\mapsto} \hat{h}(m)\otimes a$$
and
$$\oplus_{g\in
G}(a_{g})_{g}\stackrel{\hat{h}}{\mapsto} \oplus_{g\in G}(ha_{g})_{gh^{-1}},$$
respectively, and $h\circ\hat{h}=\hat{h}\circ h$ then we have that $G\times G$
acts on both sheaves and
again $\phi$ is $G\times G$-invariant.

On the other hand, let consider in $\mo_{Y}[G]$ the actions given by
$$(a_{g})_{g}\stackrel{h}{\mapsto}(a_{g})_{hg}$$
and
$$(a_{g})_{g}\stackrel{\hat{h}}{\mapsto}(a_{g})_{gh^{-1}}$$
then $h\circ\hat{h}=\hat{h}\circ h$ and again we have a $G\times G$ action,
and the isotypical decomposition of $\mo_{Y}[G]$ is given by
$$\mo_{Y}[G]=\oplus_{V\in \m I}\mo (V)\otimes \mo (V^{\vee})$$
where $\m I $ is the set of irreducible representations of $G$ over $k$.

\begin{Lema} The natural isomorphism
$$j:\mo_{Y}[G]\otimes \pib \mf\ra \oplus_{g\in G}\pib\mf$$
is $G\times G$ invariant
\end{Lema}
\textit{proof:} Let $\delta _{h}(g_{0})$ be 1 if $h=g_{0}$ and zero in other
case. Then is sufficient to proof that
$$j [h\times \hat{h}((\delta_{g}(g_{0}))_{g}\otimes m )]= h\times
\hat{h}[j((\delta_{g}(g_{0}))_{g}\otimes m )]$$
but this is immediately by direct calculation. $\diamondsuit$

\vspace{3mm}

Let $\phi_{i}$ be the map defined in \ref{paso1} Then
\begin{Prop} $\phi_{i}$
is $G$-invariant
with the right action \end{Prop}
\textit{Proof:} Recall that the morphism
$$\phi: \pib\mf\otimes\pib\mo_{X} \ra \mo_{Y}[G]\otimes \pib \mf $$
is $G\times G$ invariant. In particular we have
$$\phi: \pib\mf\otimes \mo(V)\otimes \me_{V} \ra \mo (V)\otimes \mo
(V^{\vee})\otimes \pib \mf $$ is $G\times G$ invariant. Then
$$\phi: \pib\mf\otimes \me_{V} \ra  \mo (V^{\vee})\otimes \pib \mf $$
is $G$ invariant under the right action. $\diamondsuit$

\begin{Prop}\label{paso2}
Let $X\st{\pi} Y$ be a Galois cover with Galois group $G$, and $\mf$ be a
locally free sheaf on $X$, $G$-stable, let
   $$\pib \mf= \mf_{0} \oplus \mo (V_{1})\otimes \mf_{1}\oplus
...\oplus \mo (V_{r})\otimes \mf_{r}$$
be the isotypical decomposition of $\pib \mf$. Then there exist natural
injective morphisms
 $$\mf_{i}\otimes
\me_{V_{j}}\ra \bigoplus_{l=0}^{r}\mf_{l}^{n_{l}}$$
where $$V_{i}\otimes V_{j}=\bigoplus_{l=0}^{r}V_{l}^{n_{l}}$$
\end{Prop}
\textit{Proof:} From the above proposition we have the next $G$ invariant
morphism
$$0\ra \pib
\m{F}\otimes_{\mo_{Y}}\me_{V_{j}}\st{\phi_{j}} \mo (V_{j}^{\vee}
\otimes \pib \m{F}$$

And this decompose in injective morphism
$$ (\mo(V_{i})\otimes \mf_{V_{i}})\otimes \me_{V_{j}}\ra  (\mo
(V_{j}^{\vee})\otimes \pib \m{F})$$
$$ \simeq \mo (V_{j}^{\vee})\otimes [\mf_{0} \oplus \mo
(V_{1})\otimes \mf_{1}\oplus ...\mo (V_{r})\otimes \mf_{r}]$$
$$ \simeq [\mo (V_{j}^{\vee})\otimes\mf_{0} \oplus
 \mo (V_{j}^{\vee})\otimes \mo (V_{1})\otimes \mf_{1}\oplus ... \mo
(V_{j}^{\vee}) \otimes\mo (V_{r})\otimes \mf_{r}]$$
$$ \simeq [\mo (V_{j}^{\vee})\otimes\mf_{0} \oplus
 \mo (V_{j}^{\vee}\otimes V_{1})\otimes \mf_{1}\oplus ... \mo
(V_{j}^{\vee} \otimes V_{r})\otimes \mf_{r}]$$
each one factorizing by the $V_{i}$  part. Now we introduce some notation
Let $W$ be a representation of $G$ over $k$ and
$V$ an irreducible representation; then we define $<W,V>:=$dimension that $V$
appear on $W$.

Let be $s_{l}^{i}=<V_{j}^{\vee}\otimes V_{l},V_{i}>$
then the $V_{i}$ part is determinated by
$$\mo(V_{i})\otimes [\bigoplus_{l=1}^{r}\mf_{l}^{s_{l}^{i}}]$$
from this we have the injective morphism
$$\mo(V_{i})\otimes (\mf_{i}\otimes
\me_{V_{j}})\ra
\mo(V_{i})\otimes [\bigoplus_{l=0}^{r}\mf_{l}^{s_{l}^{i}}]$$
and using the proposition 4 we have the desired injective
morphism
$$\mf_{i}\otimes
\me_{V_{j}}\ra \bigoplus_{l=0}^{r}\mf_{l}^{s_{l}^{i}}$$
on the other hand, by the theory of representation over a field we have
$<V_{j}^{\vee}\otimes V_{l},V_{i}>=
<V_{j}\otimes V_{i},V_{l}>=n_{l}$ where $$V_{j}\otimes
V_{i}=\bigoplus_{l=1}^{r}V_{l}^{n_{l}}$$
and then we conclude. $\diamondsuit$

\begin{Teo}\label{paso3} Let $X\st{\pi} Y$ be a Galois cover
with Galois group $G$, let
$$\pib \mo _{X}=\mo_{Y}\oplus\mo (V_{1})\otimes \me_{V_{1}}\oplus...\oplus\mo
(V_{r})\otimes \me_{V_{r}}$$
be the isotypical decomposition of $\pib \mo_{X}$. Then the algebra structure
determine the natural injective morphisms
 $$\me_{V_{i}}\otimes \me_{V_{j}}\ra
\bigoplus_{l=0}^{r}\me_{V_{l}}^{n_{l}} $$ where
$$V_{i}\otimes V_{j}=\bigoplus_{l=0}^{r}V_{l}^{n_{l}}.$$
Further, this are  isomorphisms in the unramified case.
\end{Teo}

For the end, we study the stability of the sheaf $\pib\mo_{X}$ in the
case that $\pi:X\ra Y$ is unramified and $X$ is a smooth variety. Then in this
case
we have that $Y$ and $\pi$ are smooth. Moreover $\pit\pib\mo_{X}\simeq
\mo_{X}^{\mid
G \mid }$ and $\pit \Omega_{Y/k}\simeq \Omega_{K/k}$ so
$$\omega_{X}=\Lambda^{dim\,
X}\Omega_{X/k}\simeq\Lambda^{dim\, X}\pit\Omega_{Y/k} \simeq \pit\Lambda^{dim\,
Y}\Omega_{Y/k} \simeq \pit \omega_{Y}.$$

Let $\mo_{Y}(1)$ be an ample line bundle on $Y$ then $\pit\mo_{Y}(1)$ is an
ample
line bundle on $X$, then we have the next

\begin{Lema} Let $\m{F}$ be an $n-$ dimensional coherent sheaf on $X$. Then
$\m{F}$
is $\mu$-polystable if an only if $\pit \m{F}$ is $\mu$-polystable.
\end{Lema}
\textit{Proof:} See \cite{len} pages 62-63.$\diamondsuit$

\begin{Teo} Let $\pi:X\ra Y$ as above. Then in the isotipycal decomposition
$$\pib\mo_{X}=\bigoplus_{V\in \m{I}}\mo(V)\otimes \m{E}_{V}$$
each $\me_{V}$ is $\mu$-stable, $\me_{V_{i}}\sim \hspace{-4mm}/\, \,
\me_{V_{j}}$ if
$i\neq j$ and $\me_{V^{*}}\simeq \me_{V}^{*}$. \end{Teo}
\textit{Proof:} By the above lemma, we have that $\pib\mo_{X}$ is
$\mu$-polistable .
Then, just we need to prove that each $\me_{V}$ is simple i.e.
$End(\me_{V_{i}})=k$,
and then $\mu$-stable. For that,
$$Hom(\pib \mo_{X},\pib \mo_{X}) \supseteq \bigoplus_{i=0}^{r}
Hom(\mo (V_{i})\otimes \me_{V_{i}},\mo (V_{i})\otimes
\me_{V_{i}})$$
$$=\bigoplus_{i=0}^{r} Hom(\me_{V_{i}},\me_{V_{i}})^{(dim\,
V_{i})^{2}}$$ now $h^{0}(\m{H}\mathit{om}(
\me_{V_{i}},\me_{V_{i}}))\geq 1$ so we have that $$h^{0}(\m{H}\mathit{om}(\pib
\mo_{X},\pib \mo_{X}))\geq \sum_{i=0}^{r}(dim\, V_{i})^{2}=\mid G\mid$$  on
the other
hand
 $$h^{0}(\pib \mo_{X}\otimes \pib
\mo_{X})=h^{0}(\pit \pib \mo_{X})= h^{0}(\mo_{X}^{\oplus \mid G\mid})=\mid
G\mid$$  and  $(\pib \mo_{X})^{\vee}\sim \pib \mo_{X}$ then
$$dim \, Hom(\pib \mo_{X},\pib \mo_{X})=\mid G\mid,$$ but this equality imply
$h^{0}(Hom(\me_{V_{i}},\me_{V_{i}}))=1$ so each $\me_{V_{i}}$ is
simple and from here is $\mu$-stable. Now observe that  $h^{0}
(Hom(\me_{V_{i}},\me_{V_{j}}))=0$ for $i\neq j$, this imply that
$\me_{V_{i}}\sim \hspace{-4mm}/\, \, \me_{V_{j}}$ if $i\neq j$.
Now just we need to see that
$(\me_{V_{i}})^{\vee}\sim \me_{V_{i}^{\vee}}$,
but this is consequence from the fact that the trace map
$\pib \mo_{X} \st{Tr} (\pib \mo_{X})^{\vee}$ is a $G$-invariant
isomorphism.$\diamondsuit$

\vspace{5mm}

\begin{flushleft}
Instituto de F\qi sica y Matem\qa ticas\\
Universidad Michoacana de\\
San Nicol\'as de Hidalgo\\
A.P. 2-82 CP 58040 \\
Morelia, Mich. M\qe xico\\
e-mail:armando@itzel.ifm.umich.mx
\end{flushleft}

\end{document}